\documentclass[11pt]{article}
\usepackage{amsmath}
\usepackage{amsmath,amssymb,amsfonts,amsthm,fancyhdr}
\usepackage{epsfig,graphicx,picins,picinpar,subfigure}
\usepackage{pstricks}
\usepackage{fancyvrb}
\usepackage[numbers,sort&compress]{natbib}
\begin{document}

\title{\textbf{Dickson's conjecture on $Z^n$\\---\small \it An equivalent form of Green-Tao's conjecture}}
\author{Shaohua Zhang}
\date{{\small  \emph{School of Mathematics, Shandong University,
Jinan,  Shandong, 250100, PRC}}}

\footnotetext [1]{E-mail address: shaohuazhang@mail.sdu.edu.cn}
\footnotetext [2]{This work was partially supported by the National
Basic Research Program (973) of China (No. 2007CB807902) and the
Natural Science Foundation of Shandong Province (No. Y2008G23).}
\maketitle

\vspace{3mm}\textbf{Abstract}

\vspace{3mm}In [1], we give Dickson's conjecture on $N^n$. In this
paper, we further give Dickson's conjecture on $Z^n$ and obtain an
equivalent form of Green-Tao's conjecture [2]. Based on our work, it
is possible to establish a general theory that several multivariable
integral polynomials on $Z^n$ represent simultaneously prime numbers
for infinitely many integral points and generalize the analogy of
Chinese Remainder Theorem in [3].

\vspace{3mm}\textbf{R\'{e}sum\'{e}}

\vspace{3mm}\textbf{Conjecture de Dickson sur $Z^n$--- Une forme
\'{e}quivalente de conjecture de Green-Tao.} Dans [1],  nous donnons
la conjecture de Dickson sur $N^n$. Dans ce document, en outre nous
accordons une conjecture de Dickson sur $Z^n$ et obtenons une forme
\'{e}quivalent de conjecture de Green-Tao [2]. Sur la base de nos
travaux, il est possible d'\'{e}tablir  une th\'{e}orie
g\'{e}n\'{e}rale que plusieurs polyn\^{o}mes int\'{e}graux
multivariables sur $Z^n$ repr\'{e}sentent simultan\'{e}ment les
nombres premiers pour un nombre infini de points entiers et de
g\'{e}n\'{e}raliser les l'analogie de Th\'{e}or\`{e}me des Restes
Chinois dans [3].

\vspace{3mm}\textbf{Keywords:} Dickson's conjecture, affine-linear
form,  Generalized Chinese Remainder Theorem, prime,  multivariable
integral polynomials, admissible map, strongly admissible map,
polynomial map on $Z^n$

\vspace{3mm}\textbf{2000 MR  Subject Classification:} 11A41, 11A99

\section{Some basic notations}
\setcounter{section}{1}\setcounter{equation}{0}

Let $Z$ be the set of integers. Denote the set of all positive prime
numbers by $P$ (we do not consider negative primes.) and denote the
set of all natural numbers or positive integers by $N$. We define
affine $n$-space over $Z$ (resp. $N$), denoted $Z^n$ (resp. $N^n$),
to be the set of all $n$-tuples of elements of $Z$ (resp. $N$). An
element $x=(x_1,...,x_n)$ in $Z^n$ is called an integral point and
$x_i$ is called the coordinates of $x$. Let $Z[x_1,...,x_n]$ be the
polynomial ring in $n$ variables over $Z$. Let $Z_n^*=\{x\in N|1\leq
x\leq n, \gcd(x,n)=1\}$ be the set of positive integers less than or
equal to $n$ that are coprime to $n$. Let $Z_n=\{x\in Z|0\leq x\leq
n-1\}$.

\section{Some basic definitions}
Let's consider the map $F: Z^n\rightarrow Z^m$ for all integral
points $x=(x_1,...,x_n)\in Z^n$, $F(x)=(f_1(x),...,f_m(x))$ for
distinct polynomials $f_1,...,f_m\in Z[x_1,...,x_n]$, where $m,n\in
N$. In this case, we call $F$ a polynomial map on $Z^n$. Call $F$ a
polynomial map on $N^n$ if $F: N^n\rightarrow Z^m$ for all integral
points $x\in N^n$, $F(x)=(f_1(x),...,f_m(x))$.  We call $F$ on $Z^n$
(resp. $N^n$ ) linear or a system of non-constant affine-linear
forms if for any $1\leq i\leq m$, $f_i(x_1,...,x_n)$ has degree 1.

\vspace{3mm}We call the polynomial map $F$ on $Z^n$ ( resp. $N^n$ )
admissible if for every positive integer $r$ there exists an
integral point $x=(x_1,...,x_n)\in Z^n$ ( resp. $N^n$ ) such that
$r$ is coprime to $f_1(x)\times...\times
f_m(x)=\prod_{i=1}^{i=m}f_i(x)$, moreover, $f_i(x)>1$ for $1\leq
i\leq m$ (namely, $F(x)\in(N\setminus \{1\})^m$). Generally
speaking, when people talk about the definition of "admissible",
they do not consider the condition "$f_i(x)>1$ for $1\leq i\leq m$".
Nevertheless, from this paper, one will see that this is very
necessary. Sometimes, in the definition of "admissible", people
regard $r$ as a prime. But, this is not always true. For example,
let's consider one-variable function $f(x)=-x^2+6$. It is easy to
test that for any prime $p$, there is an $x$ such that
$\gcd(f(x),p)=1$ and $f(x)>1$. But $f(x)$ is not "admissible".

\vspace{3mm}We call the polynomial map $F$ on $Z^n$ ( resp. $N^n$ )
strongly admissible if there is a positive integer $C$ such that for
every positive integer $k\geq C$, there exists an integral point
$x=(x_1,...,x_n)$ such that $f_1(x)>1,...,f_m(x)>1$ are all in
$Z_k^*$. We call the least positive integer $C$ such that $F$ is
strongly admissible a strongly admissible constant.

\vspace{3mm}Let $S$ be the set of all solutions of the simultaneous
equations: $$f_1(x_1,...,x_n)\in P,...,f_m(x_1,...,x_n)\in P.$$

Let $H=F(S)\in P^m$ be the image of $S$ under $F$. An element of $H$
is called a prime point. An important question is the following: Let
$F$ be a given polynomial map on $Z^n$ (resp. $N^n$), how to
determine whether $H$ is an infinite set or not?

\vspace{3mm}We say that several multivariable integral polynomials
$f_1(x),...,f_m(x)$ on $Z^n$ (resp. $N^n$) represent simultaneously
prime numbers for infinitely many integral points $x$, if for any
$1\leq i\leq m$, $f_i(x)$ itself can represent primes for infinitely
many integral points $x$, moreover, there is an infinite sequence of
integral points $(x_{11},...,x_{n1})$, ..., $(x_{1i},...,x_{ni})$,
... such that for any positive integer $r$,
$f_1(x_{1r},...,x_{nr})$,..., $f_m(x_{1r},...,x_{nr})$ represent
simultaneously primes, and for any $i\neq j$,
$f_1(x_{1i},...,x_{ni})\neq f_1(x_{1j},...,x_{nj})$, ...,
$f_m(x_{1i},...,x_{ni})\neq f_m(x_{1j},...,x_{nj})$ hold
simultaneously. In this case, we also say that the polynomial map
$F$ on $Z^n$ ( resp. $N^n$ ) represents infinitely many prime
points. One of our aims is to study the conditions that the
polynomial map $F$ on $Z^n$ represents infinitely many prime points.
In the case that $F$ is on $N^n$, we have initially considered this
problem, for the details, see [3].

\section{Introductions}
In [1], we have considered the case that $F$ on $N$ is linear and
obtained an equivalent form of Dickson's conjecture which gives the
necessary and sufficient condition that a system of non-constant
affine-linear forms on $N$ represents infinitely many prime points.
Moreover, we generalized Dickson's conjecture to the case that the
linear polynomial map $F$ is on $N^n$. In this paper,  we will
further consider the case that the linear polynomial map $F$ is on
$Z^n$ and obtain an equivalent form of Green-Tao's conjecture [2].

\vspace{3mm} Now, let's start with the work of Green and Tao. In
2006 [2], they are the first to  generalize  Dickson's conjecture
and consider the case that $F$ is a system of non-constant
affine-linear forms on $Z^n$. They noticed that if the linear
polynomial map $F=(f_1(x),...,f_m(x))$ on $Z^n$ satisfies the
conditions that for any $1\leq i\neq j\leq m$, $f_i(x)$ and $f_j(x)$
are not rational multiples of each other, then $F$ represents
infinitely many prime points. They further gave asymptotic formulae
and generalized Hardy-Littlewood Conjecture. But, the conditions
that for any $1\leq i\neq j\leq m$, $f_i(x)$ and $f_j(x)$ are not
rational multiples of each other is not sufficient and it is only
necessary. For example, $F=(f_1(x),f_2(x),f_3(x))$ never represents
prime points, where
$f_1(x)=-2x_1+3x_2-6,f_2(x)=3x_1-4x_2-6,f_3(x)=-9x_1+10x_2-6$. Based
on the idea in [3], we fix this condition as follows:

\vspace{3mm}\noindent {\bf Condition A:~~}%
For every positive integer $r$, there exists an integral point
$x=(x_1,x_2,...,x_n)$ such that for any $1\leq i\leq m$, $f_i(x)>1$
and $\gcd (f_i(x),r)=1$.

\vspace{3mm}Clearly, Condition A is necessary. We hope that it also
is sufficient. Thus, we have the following conjecture:

\vspace{3mm}\noindent {\bf Conjecture 1:~~}%
If the linear polynomial map $F$ on $Z^n$ is admissible, then $F$
represents infinitely many prime points.

\vspace{3mm}In the case that the linear polynomial map $F$ is on
$N^n$, in 2009, the author [1] obtained the equivalent form of
Dickson's conjecture and conjectured that if  $F$ on $N^n$ is
strongly admissible, then $F$ represents infinitely many prime
points. Naturally, one might further conjecture the following:

\vspace{3mm}\noindent {\bf Conjecture 2:~~}%
If the linear polynomial map $F$ on $Z^n$ is strongly admissible,
then $F$ represents infinitely many prime points.

\vspace{3mm} In order to explain the rationality of Conjecture 2, we
go back to Dickson's conjecture. In 1904, Dickson conjectured: Let
$m\geq1$, $f_i(x)=a_i+b_ix$ with $a_i$ and $b_i$ integers,
$b_i\geq1$ (for $i=1,...,m$). If there does not exist any integer
$n>1$ dividing all the products $\prod_{i=1}^{i=s}f_i(k)$, for every
integer $k$, then there exist infinitely many natural numbers $m$
such that all numbers $f_1(m),...,f_s(m)$ are primes.

\vspace{3mm} Obviously, from the condition "$b_i\geq1$",  we deduce
that Dickson only considered the linear polynomial map $F$ on $N$.
Therefore, one might call this conjecture Dickson's conjecture on
$N$. Then, What does Dickson's conjecture on $Z$ look like? One
might naively guess: Let $m\geq1$, $f_i(x)=a_i+b_ix$ with $a_i$ and
$b_i$ integers for $i=1,...,m$. If there does not exist any integer
$n>1$ dividing all the products $\prod_{i=1}^{i=s}f_i(k)$, for every
integer $k$, then there exist infinitely many natural numbers $m$
such that all numbers $f_1(m),...,f_s(m)$ are primes. Unfortunately,
it is always not true. For instance, $f_1(x)=4x-5$ and
$f_2(x)=-3x+4$ satisfy the claimed condition, but they do not
represent simultaneously primes since we do not consider negative
primes.

\vspace{3mm} It is not hard to see that Dickson's conjecture on $Z$
should be the following: Let $m\geq1$, $f_i(x)=a_i+b_ix$ with $a_i$
and $b_i$ integers, either $b_i\geq1$ or $b_i\leq-1$ (for all
$i=1,...,m$). If there does not exist any integer $n>1$ dividing all
the products $\prod_{i=1}^{i=s}f_i(k)$, for every integer $k$, then
there exist infinitely many natural numbers $m$ such that all
numbers $f_1(m),...,f_s(m)$ are primes. This can be equivalently
stated that if the linear polynomial map $F$ on $Z$ is admissible,
then $F$ represents infinitely many prime points. By the method in
[1], one can further prove that if  the linear polynomial map $F$ on
$Z$ is admissible then $F$ also is strongly admissible. Thus we get
Conjecture 2. Note that in the multivariable case, the signs of
coefficients might be disordered.

\vspace{3mm}As a toy example, we believe that $f_1(x)=2x_1-3x_2$ and
$f_2(x)=-3x_1+4x_2$ represent simultaneously primes for infinitely
many integral points $x=(x_1,x_2)\in Z^2$.

\vspace{3mm}\noindent {\bf Theorem 1:~~}%
Conjecture 1 and Conjecture 2 are equivalent. \vspace{3mm}

For the proof of theorem 1, see next section. Similarly, one can
prove that if the linear polynomial map $F$ on $N^n$ is admissible
then $F$ also is strongly admissible. Thus, we complement the work
in [1]. Based on given equivalent forms, more precisely, by finding
"strongly admissible", it is possible to establish a general theory
that several multivariable integral polynomials on $Z^n$ represent
simultaneously prime numbers for infinitely many integral points and
generalize the analogy of Chinese Remainder Theorem in [3].

\section{Proof of Theorem 1}
In order to prove that Conjecture 1 and Conjecture 2 are equivalent,
it is enough to prove that the following conditions B, C are
equivalent. Let the linear polynomial map $F(x)=(f_1(x),...,f_m(x))$
be on $Z^n$.

\vspace{3mm}\noindent {\bf Condition B:~~}%
$F$ is admissible.

\vspace{3mm}\noindent {\bf Condition C:~~}%
$F$ is strongly admissible.

\vspace{3mm}\noindent {\bf Lemma 1:~~}%
Let $F=\left\{
\begin{array}{c}
f_1(x_1,...,x_n)=a_{11}x_1+...+a_{1n}x_n +b_1\\
...........................................................\\
f_m(x_1,...,x_n)=a_{m1}x_1+...+a_{mn}x_n +b_m\\
\end{array}
\right.$ be the linear polynomial map. If $F$ is admissible, then
for any positive integer $c$, there is an integral point $x$ such
that any prime divisor of $\prod_{i=1}^{i=m}f_i(x)$ is greater than
$c$, moreover $$\left\{
\begin{array}{c}
a_{11}x_1+...+a_{1n}x_n>e=\max_{1\leq j\leq m}\{|a_{j1}|+...+|a_{jn}|\}>0\\
.....................................................................................\\
a_{m1}x_1+...+a_{mn}x_n>e=\max_{1\leq j\leq m}\{|a_{j1}|+...+|a_{jn}|\}>0\\
\end{array}
\right..$$

\vspace{3mm}\noindent {\bf Proof:~~}%
The case $b_i=0$ for each $1\leq i\leq m$ is trivial. So, let's
assume that $B=\{b_1,...,b_m\}\neq \{0\}$ and consider the number
$\alpha=\prod_{x\in B,x\neq 0}x$. By the known condition, there is
an integral point $x$ such that $$\gcd
(\prod_{i=1}^{i=m}f_i(x),\prod_{p\leq 2|\alpha| ce}p)=1,$$ moreover
$f_i(x)>1$. Notice that for each $1\leq i\leq m$, any prime divisor
of  $f_i(x)$ is greater than $2|\alpha| ce$. Therefore,
$a_{i1}x_1+...+a_{in}x_n>e$ for each $1\leq i\leq m$. Since $F$ is
admissible, hence $e>0$ and Lemma 1 holds.

\vspace{3mm}\noindent {\bf Proof that Conditions B, C are
equivalent:~~}%
Clearly, Condition C $\Rightarrow$  Condition B. Next, we prove that
Condition B $\Rightarrow$  Condition C. We write $$F=\left\{
\begin{array}{c}
f_1(x_1,...,x_n)=a_{11}x_1+...+a_{1n}x_n +b_1\\
...........................................................\\
f_m(x_1,...,x_n)=a_{m1}x_1+...+a_{mn}x_n +b_m\\
\end{array}
\right. .$$

\vspace{3mm}By lemma 1 and Condition B, we choose an integer point
$y=(y_1,...,y_n)$ such that, for $1\leq i\leq m$, $f_i(y)>1$ and
$\left\{
\begin{array}{c}
a_{11}y_1+...+a_{1n}y_n>e>0\\
................................\\
a_{m1}y_1+...+a_{mn}y_n>e>0\\
\end{array}
\right.$.

\vspace{3mm}By lemma 1 and Condition B again, we can choose an
integer point $z=(z_1,...,z_n)$ such that any prime divisor of
$\prod_{i=1}^{i=m}f_i(z)$ is greater than $\prod_{i=1}^{i=m}f_i(y)$.
We also can choose an integer point $w=(w_1,...,w_n)$ such that any
prime divisor of $\prod_{i=1}^{i=m}f_i(w)$ is greater than
$1+e+m\times \prod_{i=1}^{i=m}f_i(z)\times \max_{1\leq j\leq m}
\{|a_{j1}y_1|+...+|a_{jn}y_n|\}$.

\vspace{3mm}We claim that there is a positive integer $C\leq
\prod_{i=1}^{i=m}f_i(z)f_i(w)$ such that for every positive integer
$k\geq C$, there exists an integral point $x=(x_1,...,x_n)$ such
that $f_1(x)>1,...,f_m(x)>1$ are all in $Z_k^*$. Namely, Condition B
$\Rightarrow$  Condition C.

\vspace{3mm} Obviously, if $\gcd (k,\prod_{i=1}^{i=m}f_i(z))=1$,
then we can choose $x=z$. And if $\gcd
(k,\prod_{i=1}^{i=m}f_i(w))=1$, then we can choose $x=w$. Therefore,
it is enough to consider the case of $k=pqt$, where $p$ is a prime
divisor of $\prod_{i=1}^{i=m}f_i(z)$ and $q$ is a prime divisor of
$\prod_{i=1}^{i=m}f_i(w)$, $t$ is any positive integer.

\vspace{3mm} By our assumption that $F$ is admissible, one can
choose an integral point $v=(v_1,...,v_n)$ such that $\gcd
(t,\prod_{i=1}^{i=m}f_i(v))=1$. Write $t=p^ed$ with $\gcd (p,d)=1$.
Clearly, if $d=1$, then we can choose $x=y$ such that Condition B
$\Rightarrow$  Condition C since $\gcd
(pqt,\prod_{i=1}^{i=m}f_i(y))=1$ and $f_i(y)>1$. Hence, we might
assume that $d>1$.

\vspace{3mm}Now, we consider the coordinates. If for each $1\leq
i\leq n$, $y_i\equiv v_i (\mod d)$, then, we can similarly choose
$x=y$. So, now, we only consider the case that for some $1\leq i\leq
n$, $y_i\neq v_i (\mod d)$. Without loss of generality, we assume
that $y_1\neq v_1 (\mod d)$.

\vspace{3mm} Since $\gcd (p,d)=1$ and  $y_1\neq v_1 (\mod d)$, hence
there must be a positive integer $r_1$ with $r_1<d$ such that
$y_1+pr_1\equiv v_1 (\mod d)$. Also  from $\gcd (p,d)=1$, we deduce
that there is an integer $r_i$ with $0\leq r_i<d$ such that
$y_i+pr_i\equiv v_i (\mod d)$ for each $2\leq i\leq n$.

\vspace{3mm} Let
$x^{(h)}=(y_1+pr_1+pdhy_1,y_2+pr_2+pdhy_2,...,y_n+pr_n+pdhy_n)$ with
$h\in N$. It is easy to see that $\gcd
(pt,\prod_{i=1}^{i=m}f_i(x^{(h)}))=1$. We also have that for each
$1\leq i\leq m$,
$f_i(x^{(h)})=f_i(y)+p(a_{i1}r_1+...+a_{in}r_n)+pdh(a_{i1}y_1+...+a_{in}y_n)>
1+p(a_{i1}r_1+...+a_{in}r_n)+pdhe>1+pdhe-pd(|a_{i1}|+...+|a_{in}|)>1$.

\vspace{3mm} Finally, let $h$ range over $\{1,2,...,m+1\}$ and let's
consider the following matrix:
$$M=(m_{i,j})=\left(
\begin{matrix}
f_1(x^{(1)}), & \cdots,  & f_1(x^{(m+1)})\\
\cdots, & \cdots, & \cdots \\
f_m(x^{(1)}), & \cdots, & f_m(x^{(m+1)})
\end{matrix} \right).$$

Since $q>1+e+m\times \prod_{i=1}^{i=m}f_i(z)\times \max_{1\leq j\leq
m}\{|a_{j1}y_1|+...+|a_{jn}y_n|\}$, hence there is at most a number
which can be divided $q$ in each row of the matrix $M$. But there
are $m+1$ columns in $M$. So, there must be some $j$ with $1\leq
j\leq m+1$ such that $(\prod_{i=1}^{i=m}f_i(x^{(j)}), q)=1$. Also
notice that for each $1\leq i\leq m$, $1<f_i(x^{(j)})<pqt$. This
shows immediately that Theorem 1 holds.

\section{Further considerations}
In this paper, we mainly research the problem that the linear
polynomial map $F$ on $Z^n$ represents infinitely many prime points
and try to give a necessary and sufficient condition that $F$
represents infinitely many prime points. Based on this work and also
based on the idea in [3], it is possible to further generalize the
famous H Hypothesis as follows:

\vspace{3mm}\noindent{\bf Conjecture 3 (Generalized H Hypothesis):~~}%
Let $F$ be a polynomial map on $Z^n$. If $F$ is strongly admissible,
and there exists an integral point $y=(y_1,...,y_n)$ such that
$f_1(y)\geq C,..., f_m(y)\geq C$ are all primes, then $F$ represents
infinitely many prime points, where $C$ is the strongly admissible
constant.

\vspace{3mm}\noindent {\bf Lemma 2 (Generalized Chinese Remainder Theorem):~~}%
If $\gcd (a_i,a_j)=1$ for $1\leq i\neq j\leq m$, then the Cartesian
product $(Z_{a_1}^*)^n\times ...\times (Z_{a_m}^*)^n$ (resp.
$(Z_{a_1})^n\times ...\times (Z_{a_m})^n$) is isomorphic to
$(Z_{a_1...a_m}^*)^n$ (resp. $(Z_{a_1...a_m})^n$) for any positive
integers $a_i,m,n$.

\vspace{3mm}By Lemma 2 and the idea in [3], we further generalize
the analogy of Chinese Remainder Theorem in [3].

\vspace{3mm}\noindent {\bf Generalized  analogy of Chinese Remainder Theorem:~~}%
Let $F=(f_1,...,f_m)$ be a polynomial map on $Z^n$. If $F$
represents infinitely many prime points,  and if $\gcd (a_i,a_j)=1$
for $1\leq i\neq j\leq k$, and  there exist integral point
$x^{(j)}\in Z^n$ such that $F(x^{(j)})$ is in $(Z_{a_j}^*\setminus
\{1\})^m$ for $1\leq j\leq k$, then there exists an integral point
$z$ such that $F(z)$ is in $(Z_{a_1...a_k}^*\setminus \{1\})^m$.

\section{Acknowledgements}
I am very thankful to the referee for reading the paper, and also to
my supervisor Professor Xiaoyun Wang for her suggestions. Thank the
key lab of cryptography technology and information security in
Shandong University and the Institute for Advanced Study in Tsinghua
University, for providing me with excellent conditions.

\section{References}

\vspace{3mm}\noindent[1]  Shaohua Zhang,  Notes on Dickson's
Conjecture, available at: \\http://arxiv.org/abs/0906.3850

\vspace{3mm}\noindent[2]  Ben Green, Terence Tao, Linear equations
in primes, preprint, available at:
http://arxiv.org/abs/math/0606088, Ann. of Math. (2), in press

\vspace{3mm}\noindent[3] Shaohua Zhang, On the infinitude of some
special kinds of primes, available at:
http://arxiv.org/abs/0905.1655

\clearpage
\end{document}